\documentclass[12pt]{article}

\usepackage[colorlinks,citecolor=blue,urlcolor=blue]{hyperref}

\usepackage[english]{babel}

\usepackage{geometry}
\usepackage{amsmath,amssymb,amsfonts}

\usepackage[numbers]{natbib}
\usepackage{bm}
\usepackage{dsfont}
\usepackage{mathtools}
\usepackage{graphicx}
\usepackage{tabularx}
\usepackage{booktabs}
\usepackage{xspace}

\usepackage{amsthm}
\newtheorem{theorem}{Theorem}
\newtheorem{lemma}{Lemma}

\newenvironment{IEEEproof}{\proof}{\endproof}
\newenvironment{IEEEkeywords}
    {\textbf{Keywords:}}
    {}

\newenvironment{MSC}
    {\textbf{\href{https://mathscinet.ams.org/mathscinet/msc/msc2020.html}{MSC2020}:}}
    {}

\newcommand{\IEEEQED}{\qed}

\newcommand{\IEEEauthorblockN}{}

\newcommand{\dir}{\text{Dir}}
\newcommand{\mn}{\text{MN}}

\newcommand{\unif}{\text{U}}

\newcommand{\ber}{\text{Ber}}
\newcommand{\err}{\operatorname{Err}}
\newcommand{\Ham}{\operatorname{Ham}}
\newcommand{\ace}{d_\text{ACE}}

\newcommand{\cZ}{\pazocal{Z}}
\newcommand{\cX}{\mathfrak{X}}

\newcommand{\Mir}{d_{\rm Mir}}
\newcommand{\Rand}{d_{\rm Rand}}

\newcommand{\Ber}{\ber}

\DeclarePairedDelimiter\abs{|}{|}

\newcommand{\widesim}[2][1.5]{
  \mathrel{\overset{#2}{\scalebox{#1}[1]{$\sim$}}}
}

\usepackage{calrsfs}
\DeclareMathAlphabet{\pazocal}{OMS}{zplm}{m}{n}

\makeatletter
\newcommand*{\rom}[1]{\expandafter\@slowromancap\romannumeral #1@}
\makeatother

\usepackage{enumitem,xcolor}
\usepackage{cleveref}
\crefname{task}{Task}{Tasks}
\crefname{case}{Case}{cases}

\begin{document}

\title{Consistent Bayesian community recovery in multilayer networks}

\author{\IEEEauthorblockN{Kalle Alaluusua and Lasse Leskel\"a}\\
\textit{Department of Mathematics and Systems Analysis} \\
\textit{Aalto University}\\
Espoo, Finland \\
\href{mailto:kalle.alaluusua@aalto.fi}{kalle.alaluusua@aalto.fi}; \href{mailto:lasse.leskela@aalto.fi}{lasse.leskela@aalto.fi}}

\maketitle
\thispagestyle{plain}
\pagestyle{plain}

\begin{abstract}
Revealing underlying relations between nodes in a network is one of the most important tasks in network analysis. Using tools and techniques from a variety of disciplines, many community recovery methods have been developed for different scenarios. Despite the recent interest on community recovery in multilayer networks, theoretical results on the accuracy of the estimates are few and far between. Given a multilayer, e.g. temporal, network and a multilayer stochastic block model, we derive bounds for sufficient separation between intra- and inter-block connectivity parameters to achieve posterior exact and almost exact community recovery. These conditions are comparable to a well known threshold for community detection by a single-layer stochastic block model. A simulation study shows that the derived bounds translate to classification accuracy that improves as the number of observed layers increases.

\vspace{12pt} \noindent 

\noindent \begin{MSC}
05C80, 60B10, 62F15, 62H30, 90B15, 94C15
\end{MSC}

\end{abstract}

\section{Introduction}

Data sets in many application domains, such as physics, sociology, computer science, economics, epidemiology and neuroscience, consist of pairwise interactions. An important unsupervised learning problem is to infer latent community memberships from observed pair interactions, when nodes in the same community are, in some sense, more similar to each other than to the other nodes. This task is commonly known as community recovery, community detection, or clustering.

Community recovery is typically approached by fitting a generative model, such as a stochastic block model (SBM) \cite{holland_laskey_leinhardt_1983}, on the observed network data. The stochastic block model is a probability distribution on the space of adjacency matrices, where the link probability between two nodes is solely determined by the assignment of the nodes into communities, also referred to as blocks. Thus, once the community assignment and the matrix of link probabilities is known, it is easy to sample networks from the model or evaluate the likelihood of given data.

For any given network, the performance of a community detection algorithm can be evaluated by inspecting metrics such as the classification accuracy. This leads to procedures that seem to work well in practice while no claims about their asymptotic properties are made. This is in particular the case with Bayesian methods \cite{Hofman_Wiggins_2008, Peixoto_2019} that have gained popularity due to their flexibility and easy adaptability into various modelling contexts and data types. However, without exercising proper care or judgment, assigning a prior probability distribution on a large parameter space typical for community recovery problems risks failure. Quantifying how well a particular prior combined with a large data set succeeds in outputting a posterior distribution well concentrated near the corrected parameter value, is referred to as Bayesian consistency. In the context of single-layer networks, recent studies include \cite{vanderpas_vandervaart_2018, ghosh_pati_bhattacharya_2019, kleijn_van_waaij_2021,jiang_tokdar_2021, li_zhou_tsui_wei_ji_2021}.

The stochastic block model has been expanded to model vector-valued or equivalently multilayer, e.g. temporal, interactions between nodes. SBM variants with overlapping communities include those of \cite{
herlau_morup_schmidt_2013, xu_hero_2014,han_xu_airoldi_2015, matias_miele_2017}. Dynamic community recovery combines aspects of time series, where the time increases, and machine learning, where the number of observations increases. These problems are typically characterized by multilayer data, e.g. a three-way adjacency tensor indexed by nodes and a time parameter. Despite recent interest, the number of theoretical results in this direction has so far been rather limited \cite{ghasemian_zhang_clauset_moore_peel_2016, Barucca_Lillo_Mazzarisi_Tantari_2018, Bhattacharyya_Chatterjee_2020-04-06,durante_dunson_vogelstein_2017,lei_chen_lynch_2020, Avrachenkov_Dreveton_Leskela_2022}. In particular, there exists little research on Bayesian consistency for multilayer network models.  

This paper contributes to the field by deriving bounds for sufficient separation between  intra- and inter-block connectivity parameters for consistent Bayesian community recovery in multilayer networks. The theoretical results are demonstrated on simulated networks of small to moderate size. The simulation study shows that the classification accuracy improves as more network layers are observed.

\subsection{Notation}
We denote by $[n]$ the set $\{1,2,\dots,n\}$. The notation $P f$ is an abbreviation of $\int fdP$. When $\Pi$ denotes a prior distribution, $P\Pi$ is the expected posterior probability when data are sampled from $P$. For probability measures $F,G$ with densities $f,g$ relative to measure $\mu$, define ${\rho_\alpha(f\mid\mid g)} =  \int f^\alpha g^{1-\alpha}d\mu$, $\alpha \in (0,1)$. Denote by $D_{\alpha}(f\mid \mid g) = (\alpha-1)^{-1} \log\rho_\alpha(f\mid \mid g)$ the R\'enyi divergence of order $\alpha$ between $f$ and $g$ (see, e.g., \cite{vanErven_Harremoes_2014}). We consider intra- and inter-block interaction distributions $f = \otimes_{t=1}^T\ber_{p_t}$ and $g = \otimes_{t=1}^T\ber_{q_t}$, respectively, where $\ber_p(x) = (1-p)^{1-x} p^x$ denotes the Bernoulli distribution with mean $p$. Finally, $I_T \coloneqq D_{1/2}(f,g) = \sum_{t=1}^TI(p_t,q_t)$, where we define $I(p_t,q_t) = D_{1/2}(\ber_{p_t},\ber_{q_t})$.

\subsection{Multilayer stochastic block model}
\label{sec:dynamic_sbm}

Consider a multilayer SBM with $N$ nodes, $T$ layers, $K = 2$ blocks, intra-block link probabilities $p_1,\dots, p_T$, and inter-block link probabilities $q_1,\dots, q_T$. Denote the set of block structures by
$\cZ = \{z: [N] \to [K]\}$. Let $Q^{(t)}$ be a $K\times K$ matrix such that $Q^{(t)}_{z(i)z(j)} = p_t$ for $z(i)=z(j)$ and $Q^{(t)}_{z(i)z(j)} = q_t$ otherwise. Denote the space of observations by
\begin{align*}
 \cX
 = \Big\{&X: [N] \times [N] \times [T] \to \{0,1\}:
    X^{(t)}_{ij} = X^{(t)}_{ji}, \, X^{(t)}_{ii}=0 \ \text{for all $i,j,t$} \Big\}.
\end{align*}
Given a node labelling $z$, the observation is distributed according to the probability measure $P_z$ on $\cX$ defined by
\begin{align}
\begin{split}
 P_z(X)
 &= \prod_{1 \le i < j \le N}
 \bigg(F_{z(i)z(j)}(X_{ij}) \bigg)
 = \prod_{1 \le i < j \le N}\prod_{1 \le t \le T}
 \bigg(\ber_{Q^{(t)}_{z(i)z(j)}}(X^{(t)}_{ij}) \bigg)
\end{split}
\label{eq:adjacency_dynamic}
\end{align}
where $F$ is be the matrix of intra- and inter-block interaction distributions $f = \otimes_{t=1}^T\ber(p_t)$ and $g = \otimes_{t=1}^T\ber(q_t)$ such that $F_{z(i)z(j)} = f$ for $z(i)=z(j)$ and $F_{z(i)z(j)} = g$ otherwise.

\subsection{Bayesian inference}

Given a prior probability distribution $\Pi$ on $\cZ$ and an observation $X \in \cX$, denote the corresponding posterior distribution by
\[
 \Pi_X(w)
 = \frac{\Pi(w) P_w(X)}{\sum_{w'} \Pi(w') P_{w'}(X)},
 \qquad w \in \cZ.
\]
We consider $\Pi_X$ as a Bayesian distributional estimate of an unknown block structure, from which point estimates can be derived for example by taking a mode.  The accuracy of such estimates can be analysed using a frequentist viewpoint where we assume that the observed data tensor $X$ is sampled from a model $P_z$ with true block structure $z$, and we compute the expected mass that the posterior distribution assigns near the true value according to
\begin{equation}
 \label{eq:Err}
 \err_z(r)
 = P_{z} \Pi_X\{ w: \ace(w,z) > r \}.
\end{equation}
Here
$
 \ace(w,z)
 = \min\{ \Ham(w,z), \Ham(w, \tilde z)\}
$ 
denotes the absolute classification error computed 
using the Hamming distance $\Ham(w,z) = \sum_{i=1}^N (1-\delta_{w(i)z(i)})$, when $\tilde z$ is the modification of $z$ obtained by swapping the labels 1 and 2.

\subsection{Large-scale recovery}

Large-scale data regimes can be modelled using a sequence of models in which the model parameters $(N, T, p_t, q_t)$ as well as the spaces $\cX,\cZ$, the true block structure $z$, and the distributions $\Pi$ and $\Pi_X$ all depend on a scale parameter $\nu=1,2,\dots$ which is omitted from the notation for clarity. In a large-scale regime, we say that the posterior distribution \textit{exactly recovers $z$} if the error defined in (\ref{eq:Err}) satisfies
\[
 \err_z(0) = o(1).
\]
Note that $\ace(w,z) = 0$ if and only if $w \in \{z, \tilde z\}$. Posterior exact recovery hence means that most of the posterior mass is concentrated exactly at the set $\{z,\tilde z\}$ corresponding to the true unlabelled block structure. Similarly, we say that the posterior \textit{almost exactly recovers $z$} if
\[
 \err_z(\epsilon N) = o(1)
\]
for every scale-independent constant $\epsilon > 0$. Almost exact recovery means that with high probability, most of the posterior mass is concentrated on the set of block structures $w$ for which the relative classification error $N^{-1} \ace(w,z)$ is at most $\epsilon$.

\section{Information-theoretic thresholds}
\label{sec:consistent_dynamic}

To gain understanding on how the increase in the number of network snapshots (or layers) affects the difficulty of a community recovery problem, we derive sufficient conditions for consistent community recovery, which we compare to existing literature.

\subsection{Main results}

Recovering an unknown block structure is possible only if the
link probabilities $p_t$ and $q_t$ differ sufficiently from each other. For learning from single-layer observations \cite{abbe_bandeira_hall_2015, Mossel_Neeman_Sly_2016}, a sharp information quantity for characterising recoverability is the R\'enyi divergence of order 1/2 between Bernoulli distributions $\Ber_{p_t}$ and $\Ber_{q_t}$ given by
\[
 I(p_t, q_t)
 = (1-p_t)^{1/2} (1-q_t)^{1/2} + p_t^{1/2} q_t^{1/2}.
\]
Sparse networks are often modelled by assuming that $p_t = a_t \rho$ and $q_t = b_t \rho$ for scale-independent constants $a_t \ne b_t$ and overall link density $\rho=o(1)$, for example $\rho = N^{-1}$ (constant average degree) or $\rho = \frac{\log N}{N}$ (logarithmic average degree). In such case Taylor expansions show that
\[
 I(p_t,q_t)
 = \left( \sqrt{a_t} - \sqrt{b_t} \right)^2 \rho
 + O(\rho^2).
\]

The following theorems characterise posterior recovery from multilayer network data in terms of
\begin{equation}
 \label{eq:IT}
 I_T = \sum_{t=1}^T I(p_t,q_t).
\end{equation}

\begin{theorem}
\label{the:AlmostExactRecovery}
If $I_T \gg N^{-1}$, then the posterior distribution corresponding to the uniform prior on the set of all block structures almost exactly recovers any particular block structure $z$.
\end{theorem}

For $p_t = p$ and $q_t = q$, Theorem~\ref{the:AlmostExactRecovery} shows that $I(p,q) \gg (NT)^{-1}$ suffices for almost exact recovery. Especially, we see that almost exact recovery may be achievable for a bounded number of nodes if the number of layers $T$ is large. Reference \cite{durante_dunson_vogelstein_2017} arrives to a similar conclusion in the context of latent space models. When we view $T$ as time, we see that the product $NT$ indicates that, from an information-theoretic point of view, observing one new node in the network is equally informative to observing one new time slot.

The following theorem characterises  exact recovery of posterior distributions corresponding to noninformative priors.

\begin{theorem}
\label{the:ExactRecovery}
If $I_T \ge (2+\delta) \frac{\log N}{N}$ for some scale-independent constant $\delta > 0$, then the posterior distribution corresponding to the uniform prior on the set of all block structures exactly recovers any particular block structure $z$.
\end{theorem}

We believe that the sufficient condition in Theorem~\ref{the:ExactRecovery} is sharp because for $T=1$, it is known \cite{Mossel_Neeman_Sly_2016} that exact recovery (in a frequentist sense) is impossible when 
$I_T \le (2-\delta) \frac{\log N}{N}$. A related result by \cite{lei_chen_lynch_2020} shows that observing multiple network layers allows for consistent community detection (by a least squares estimator) from a sparser network. Specialised to $T=1$,
Theorem~\ref{the:ExactRecovery} also improves the result of \cite{kleijn_van_waaij_2021} who showed that 
$I_1 \ge (4+\delta) \frac{\log N}{N}$ is sufficient for posterior exact recovery in single-layer networks. In a frequentist setting, \cite{Avrachenkov_Dreveton_Leskela_2022} presents consistency thresholds for multilayer SBMs, which are comparable to those of Theorems \ref{the:AlmostExactRecovery} and \ref{the:ExactRecovery}.

\section{Simulation study}
\label{sec:simulation_study}
We perform a simulation study to examine the effect of the number of observed network layers on the classification accuracy of a community recovery algorithm. In particular, we study the performance of an extension of the Gibbs sampler by \cite{nowicki_snijders_2001} that takes as an input a tensor of independent and identically distributed adjacency matrices. The source code for replicating these experiments is available at \href{https://github.com/kalaluusua}{github.com/kalaluusua}.

\subsection{Posterior sampler}
The dynamic SBM introduced in this section adopts the conjugate priors by \cite{nowicki_snijders_2001},
\begin{align}
\begin{split}
    &Q_{ab} \widesim[2]{\text{iid}}\unif[0,1],~ 1\leq a \leq b\leq K,\\
    &z(i) \widesim[2]{\text{iid}} \mn_K(1;\theta),~ i\in [N]\\
    &\theta \sim \dir(K;\alpha).
\end{split}
\label{prior:dynamic}
\end{align}
where $Q_{ab}$ is the link probability between the blocks $a$ and $b$, $\mn_k(n; p_1,\dots,p_k)$ is the $k$-variate multinomial distribution with $n$ trials, and $\dir(k;\alpha_1,\dots,\alpha_k)$ is the $k$-dimensional Dirichlet distribution. The associated network follows the distribution (\ref{eq:adjacency_dynamic}), where $Q^{(t)} = Q$ for all $t$.

To take advantage of the increased number of observed layers, we adjust the likelihood function accordingly, which by the independence of $X^{(t)}$ yields $\Pi(X\mid z,Q) = \prod_{t=1}^T\pazocal{L}(z, X^{(t)})$, where
\begin{alignat}{2}
    \pazocal{L}(z, X^{(t)}) = \prod_{1\leq a\leq b\leq K}Q_{ab}^{O_{ab}(z,X^{(t)})}(1-Q_{ab})^{n_{ab}(z) - O_{ab}(z,X^{(t)})}
\label{eq:likelihood_dynamic}
\end{alignat}
where $O_{ab}$ is the number of links between communities $a$ and $b$, and $n_{ab}$ is the maximum number of links that can be formed between communities $a$ and $b$.

We propose a dynamic posterior sampler that is an extension of the Gibbs sampler introduced in \cite{nowicki_snijders_2001}. The sampler approximates the posterior densities of $(\theta, Q), z(1), \dots, z(N)$, where $(\theta, Q)$ is treated as a single random vector with the prior density $\Pi(\theta, Q)$. Given the likelihood function (\ref{eq:likelihood_dynamic}), the conditional distribution of $z(i)$ becomes
\begin{alignat*}{2}
    &~\Pi(z(i) = a\mid X,P,\theta,z_{-i}) = C\theta_a\prod_{t=1}^T\pazocal{L}^{(i)}(z_{-i}, X),
\end{alignat*}
where $z_{-i} \coloneqq (z(j))_{j\neq i}$ and
\begin{align*}
     \pazocal{L}^{(i)}(z_{-i}, X) = \prod_{b=1}^KQ_{ab}^{O^{(i)}_{b}(z_{-i}, X^{(t)})}(1-Q_{ab})^{n^{(i)}_{b}(z_{-i}) - O^{(i)}_{b}(z_{-i}, X^{(t)})},
\end{align*}
such that $C$ is a constant independent of $a$, $O^{(i)}_{b}$ is the number of $i$-$v$ links such that the node $v\in G$ is assigned to community $b$, and $n^{(i)}_{b}$ is the number of nodes $v\neq i$ in community $b$. The posterior distribution $\Pi(\theta, Q\mid z,X)$ is given by independent Dirichlet distributions with parameters
\begin{alignat*}{2}
    &(n_a + \alpha_a)_{a\in [K]} && \text{ for } \theta, 
    \\
    &\left(\sum_{t = 1}^T O_{ab}(X^{(t)})+1,Tn_{ab} - \sum_{t = 1}^T O_{ab}(X^{(t)}) + 1\right) && \text{ for } Q_{ab},
\end{alignat*}
where $1\leq a \leq b \leq K$. The posterior mode of $Q_{ab}$ becomes $1/T \sum_t O_{ab}(X^{(t)})/n_{ab}$, the average proportion of $ab$-links over all layers. This is an extension of the block constant least squares estimator used extensively in literature \cite{gao_lu_zhou_2015, klopp_tsybakov_verzelen_2017, vanderpas_vandervaart_2018}.

\subsection{Simulation design}

We study the community detection performance of the dynamic SBM on synthetically generated networks of $K = 2$ communities of size $N = 100$. We choose the following cases for link probabilities:
\begin{enumerate}[leftmargin=2cm, label=Case \arabic*:, ref=\arabic*]
\crefalias{enumi}{case}
    \item $p = 0.3$ and $q = 0.2$; \label{one}
    \item $p = 0.15$ and $q = 0.1$, \label{two}
\end{enumerate}
where $p$ and $q$ are the intra- and inter-block link probabilities, respectively. In both cases we observe 10 synthetically generated networks with $T \in \{1, 3, 5, 7\}$ independent and identically distributed network layers. To control the sources of variation, the 10 synthetically generated networks share a community structure $z_0$ where the nodes are deterministically and uniformly assigned into $K$ communities. Finally, the networks are generated from (\ref{eq:adjacency_dynamic}) with $Q_0 = \left[\begin{smallmatrix}
 p & q \\ q & p
\end{smallmatrix}\right]$.

To construct a single point estimate of community assignment from the Markov chain generated by the Gibbs sampler, we employ a method presented in \cite{dahl_2006}. The method uses information from all the community assignments and selects a certain average assignment. In practice, we average over the last $100$ members of the sequence of $1100$ states, and allow for an initial burn-in period of $1000$ initial iterations before stationary is reached.

To evaluate the accuracy of our estimate given the underlying community structure, we use the Hubert-Arabie adjusted Rand index \cite{rand_1971,hubert_arabie_1985}, which is a measure of similarity between two community assignments. The index is one when the assignments are identical and zero when they are independent. Since the definition disregards the relative community sizes, it tends to represent the level of agreement among large communities. However, in our experiment the average community sizes are close to being equal, and use of the index is justified.

\subsection{Simulation results}
Table \ref{tab:results} depicts the averages and the standard deviations of the classification error $\ace/N$ and the adjusted Rand index of our estimate, given the number of observed network layers $T$ and the link probabilities $p,q$. For each pair $(T, \text{Case } i)$, the 10 simulated networks have $N = 100$ nodes with $K = 2$ communities, and connectivity parameters that vary between \cref{one} and \cref{two}. The priors are described in (\ref{prior:dynamic}). The number of communities is known and $\alpha = (200,200).$ From now on will refer to the average adjusted Rand index by accuracy. The table also depicts the standard deviation of each associated vector of estimates. \cref{two} is more difficult than \cref{one} since it invokes sparser networks where inter-block link probability is very close to the intra-block link probability. The R\'enyi divergence $I$ values corresponding to \cref{one,two} are 0.010 and 0.006, respectively. We expect that the community detection performance improves as $T$ increases and is overall worse in \cref{two}. This is precisely what Table \ref{tab:results} shows.

\begin{table}[!htb]
\centering
\caption{Average classification error (CE) and average adjusted Rand index (AR) of the community assignment estimate with the corresponding standard deviation in parenthesis.}
\begin{tabular}{ccccccc}\toprule
& \multicolumn{2}{c}{\cref{one}} & \multicolumn{2}{c}{\cref{two}}
\\\cmidrule(lr){2-3}\cmidrule(lr){4-5}
$T$ & CE  & AR & CE & AR \\\midrule
1 & 0.34$_{(0.13)}$ & 0.16$_{(0.20)}$ & 0.45$_{(0.03)}$ & 0.01$_{(0.01)}$ \\
3 & 0.03$_{(0.01)}$ & 0.90$_{(0.05)}$ & 0.30$_{(0.15)}$ & 0.23$_{(0.23)}$ \\
5 & 0.01$_{(0.01)}$ & 0.98$_{(0.03)}$ & 0.09$_{(0.09)}$ & 0.70$_{(0.23)}$ \\
7 & 0.00$_{(0.01)}$ & 0.99$_{(0.03)}$ & 0.04$_{(0.02)}$ & 0.85$_{(0.08)}$ \\\bottomrule
\end{tabular}
\label{tab:results}
\end{table}

When we inspect each case individually, we observe that the classification error decreases and the adjusted Rand index increases as $T$ increases. When $T=1$, the sampler misclassifies on average third of the observations in \cref{one} (with an accuracy of 0.16) and nearly half of the observations (with an accuracy of 0.01) in \cref{two}. Recall that due to symmetry the domain of the classification error is $[0,1/2]$ and a Rand index of 0.01 implies that the assignments are effectively independent. The relatively large standard deviations in \cref{one} imply that, despite the bad overall performance, given a favourable initial values the sampler may correctly classify a large proportion of the nodes. Relatively small standard deviations in \cref{two} imply that this is unlikely when the problem is more difficult. In \cref{one} the accuracy improves rapidly as the number of observed layers increases; when $T=3$, the sampler is very likely to classify all but few nodes correctly, while $T$ values of 5 and 7 lead to near perfect accuracy. In \cref{two} the increase in accuracy is more muted but nevertheless apparent; when $T=3$, the accuracy of the sampler resembles that in \cref{one} with $T=1$, and when $T=7$, it resembles that in \cref{one} with $T=3$.

\section{Final remarks}
There are many important questions left unanswered, including the question of whether the sufficient conditions presented in Section \ref{sec:consistent_dynamic} are also necessary. Moreover, our analysis is limited to stochastic block models with two communities. Generalizing the results for $K > 2$ communities, with $K$ possibly unknown, is left for future work. Another interesting research direction would be to further examine community detection in networks that vary over time.

\section{Proofs}

\subsection{Mirkin distance}

The Mirkin distance between $z,w: [N] \to [K]$ is defined by
\begin{equation}
 \label{eq:Mir}
 \Mir(z,w)
 = 2(M_{01}+M_{10}),
\end{equation}
where $M_{ab}$ is the number of unordered pairs $\{i,j\}$ such that $\delta_{z(i) z(j)} = a$ and $\delta_{w(i) w(j)} = b$.  The Mirkin distance is one of the common pair-counting based cluster validity indices \cite{Gosgens_Tikhonov_Prokhorenkova_2021, Lei_etal_2017}, and it is related to the Rand index by $\Mir(z,w) = N(N-1) (1-\Rand(z, w))$.

\begin{lemma}
\label{the:MirkinLB}
For $K=2$, the Mirkin distance $M = \Mir(z,w)$, the Hamming distance $H = \Ham(z,w)$, and the absolute classification error $A = \ace(z,w)$ are related by
\begin{equation}
 \label{eq:MirkinHamming}
 M
 = 2 (N - H) H
 = 2 (N - A) A.
\end{equation}
\end{lemma}

\begin{IEEEproof}
Denote the blocks under focus by $C_k = \{i: z(i) = k\}$ and $C'_\ell = \{i: w(i) = \ell\}$. Also denote $N_{k\ell} = \abs{C_{k} \cap C'_\ell}$. For $K=2$, we find that
\begin{align*}
 M_{01} &= N_{11}N_{21} + N_{12} N_{22}, \\
 M_{10} &= N_{11}N_{12} + N_{21} N_{22}.
\end{align*}
By summing these, we find that
\[
 \frac12 \Mir(z,w)
 = M_{01} + M_{10}
 = (N_{11} + N_{22})( N_{12} + N_{21}).
\]
The first equality in \eqref{eq:MirkinHamming} follows by noting that $H = N_{12} + N_{21}$ and $N_{11} + N_{22} = N-H$. The second equality in \eqref{eq:MirkinHamming} follows by noting  that for $K=2$, the group of permutations only contains the identity map and the transposition $\tau_{12}$ which swaps 1 and 2.  In this case we find that $\Ham(z, \tau\circ w) = N - \Ham(z, w)$. Therefore, the term $N(N-H)$ in \eqref{eq:MirkinHamming} remains invariant if we replace $w$ by $\tau \circ w$.
\end{IEEEproof}

\subsection{Upper bound on posterior mass}

The following is a generalised version of \cite[Proposition 3.1:(ii)]{kleijn_van_waaij_2021} for multilayer SBMs.

\begin{lemma}
\label{the:PosteriorUB}
For any $z \in \cZ$ and any $S \subset \cZ$ not containing $z$, the expected posterior mass relative to a prior distribution $\Pi$ on $\cZ$ is bounded by
\begin{equation}
 P_z \Pi_X(S)
 \leq \sum_{w \in S} \left( \frac{\Pi(w)}{\Pi(z)} \right)^{1/2} e^{-\frac14 I_T \Mir(z,w)},
 \label{eq:posterior_contraction}
\end{equation}
where $I_T$ is defined by \eqref{eq:IT} and  $\Mir(z,w)$ by \eqref{eq:Mir}.
\end{lemma}
\begin{IEEEproof}
By \cite[Proposition D.1]{Ghosal_VanDerVaart_2017}, it follows that for any $z$ and $w$, the likelihood ratio test $\phi_{zw}(X) = \mathds{1}( \frac{P_w(X)}{P_z(X)} > \frac{\Pi(z)}{\Pi(w)} )$ satisfies
\[
 \Pi(z) P_z \phi_{zw} + \Pi(w) P_w(1-\phi_{zw}) 
 \leq \frac{\Pi(z)^\alpha}{\Pi(w)^{\alpha-1}} \rho_\alpha(P_z \mid\mid P_w)
\]
for all $0 < \alpha < 1$. By dividing both sides by $\Pi(z)$, it follows that
\begin{equation}
 \label{eq:KeyPosteriorBound}
 P_z \phi_{zw} + \frac{\Pi(w)}{\Pi(z)} P_w(1-\phi_{zw}) 
 \leq \left( \frac{\Pi(w)}{\Pi(z)} \right)^{1-\alpha} \rho_\alpha(P_z \mid\mid P_w)
\end{equation}
Let us define $\phi_z(X) = \max_{w \in S} \phi_{zw}(X)$.  
By \cite[Lemma 2.2]{Kleijn_2021}, we have
\begin{align*}
 P_z \Pi_X(S)
 \leq P_z \phi_z + \sum_{w \in S} \frac{\Pi(w)}{\Pi(z)} P_w(1-\phi_z).
\end{align*}
Then $P_z \phi_z \le \sum_{w \in S} P_z \phi_{zw}$ and $1-\phi_z \le 1-\phi_{zw}$ for all $w \in S$, so it follows by \eqref{eq:KeyPosteriorBound} that
\begin{align}
\begin{split}
 P_z \Pi_X(S)
 &\leq \sum_{w \in S} P_z \phi_{zw} + \sum_{w \in S} \frac{\Pi(w)}{\Pi(z)} P_w(1-\phi_{zw}) \\
 &= \sum_{w \in S} \left( P_z \phi_{zw} + \frac{\Pi(w)}{\Pi(z)} P_w(1-\phi_{zw}) \right) \\
 &\leq \sum_{w \in S} \left( \frac{\Pi(w)}{\Pi(z)} \right)^{1-\alpha} \rho_\alpha(P_z \mid\mid P_w).
 \end{split}
 \label{eq:posterior_contraction_intermed}
\end{align}

Recall that by (\ref{eq:adjacency_dynamic}), $P_z = \prod_{1 \le i < j \le N} F_{z(i)z(j)}$. Because R\'enyi divergence is linear with respect to products, it follows that
$
 D_\alpha(P_z \mid\mid P_w)
 = \sum_{1 \le i < j \le N} D_\alpha( F_{z(i)z(j)} \mid\mid F_{w(i)w(j)} ).
$
By definition of $F$, it now follows that
\[
 D_\alpha(P_z \mid\mid P_w)
 = M_{01} D_\alpha( g \mid\mid f ) + M_{10} D_\alpha( f \mid\mid g ).
\]
By setting $\alpha=\frac12$ and recalling definition \eqref{eq:Mir}, we find that
\begin{equation}
 D_{1/2}(P_z \mid\mid P_w)
 = \frac12 \Mir(z,w) I_T.
 \label{eq:mirkin_renyi}
\end{equation}
Inequality (\ref{eq:posterior_contraction}) follows by combining (\ref{eq:posterior_contraction_intermed}) and (\ref{eq:mirkin_renyi}).
\end{IEEEproof}

\subsection{Preliminary estimates}

 Denote $B_{z,r} = \{w: \ace(z,w) \le r\}$ and $S_{z,k} = \{w: \ace(z,w) = k\}$. By combining Lemma~\ref{the:PosteriorUB} and Lemma~\ref{the:MirkinLB}, and assuming that $\Pi$ is uniform\footnote{
It suffices to assume that $\Pi$ restricted to $S$ is uniform, and we might relax this assumption rather easily.
},
it follows that error quantity defined by \eqref{eq:Err} is bounded by
\begin{align*}
 \err_z(r)
 &= P_z \Pi_X(B_{z,r}^c) \\
 &\leq \sum_{w \in B_{z,r}^c} \left( \frac{\Pi(w)}{\Pi(z)} \right)^{1/2} e^{-\frac14 \Mir(z,w) I_T} \\
 &= \sum_{w \in B_{z,r}^c} e^{- \frac12 I_T (N-\ace(z,w)) \ace(z,w)} \\
 &= \sum_{r < k \le N/2} \abs{S_{z,k} } \, e^{- \frac12 I_T (N-k) k }.
\end{align*}
By noting that $\abs{S_{z,k}} \le 2 \binom{N}{k}$, it follows that
\begin{equation}
 \label{eq:BallComplement}
 \err_z(r)
 \leq 2 \sum_{r < k \le N/2} \binom{N}{k} \, e^{-\frac12 I_T (N-k) k}.
\end{equation}

\subsection{Proof of Theorem~\ref{the:AlmostExactRecovery}}

Recall that%
\footnote{
Because $\frac{k^k}{k!} \le \sum_{s=0}^\infty \frac{k^s}{s!} = e^k$, we see that $\binom{N}{k} \le \frac{N^k}{k!} \le (\frac{eN}{k})^k$.
}
$\abs{S_{z,k}} \le 2 \binom{N}{k} \le 2 (\frac{eN}{k})^k \le 2 (\frac{eN}{r})^k$ for all $r < k \le N/2$.
Together with the bound $(N-k) \ge N/2$, we find that by applying \eqref{eq:BallComplement} that
\begin{align*}
 \err_z(r)
 \leq 2 \sum_{r < k \le N/2} \left( \frac{eN}{r} \right)^k \, e^{-\frac14 N I_T k}
 \leq 2 \sum_{k=1}^\infty b_r^k,
 \end{align*}
where $b_r = \frac{eN}{r} e^{-\frac14 N I_T}$. For $r = \epsilon N$ for $\epsilon > 0$ being a scale-independent constant, we see that
$
 b_r
 \to 0
$ due to $I_T \gg N^{-1}$. In light of the above inequality, it follows that $ \err_z(\epsilon N) \to 0$ for every scale-independent constant $\epsilon > 0$.  Hence Theorem~\ref{the:AlmostExactRecovery} is valid. \hfill
\IEEEQED

\subsection{Proof of Theorem~\ref{the:ExactRecovery}}

To prove Theorem~\ref{the:ExactRecovery}, we will conduct a more careful analysis by splitting the sum in \eqref{eq:BallComplement} at $\ell = N^{2/3}$.

 For $1 \le k \le \ell$, we apply the inequalities $N-k \ge N-\ell$ and $\binom{N}{k} \le \frac{N^k}{k!}$ to conclude that
$
 \binom{N}{k} e^{-\frac12 (N-k) k I_T}
 \le \frac{a_\ell^k}{k!},
$
where $a_\ell = N e^{-\frac12 (N-\ell) I_T}$. For $\ell < k \le N/2$, we apply the same bounds as in the proof of Theorem~\ref{the:AlmostExactRecovery}, to conclude that
$\binom{N}{k} e^{-\frac12 (N-k) k I} \le b_\ell^k$
where $b_\ell = (\frac{eN}{\ell}) e^{-\frac14 N I_T}$. It follows by applying \eqref{eq:BallComplement} with $r=0$ that
\begin{align*}
 \err_z(0)
 &\leq 2 \sum_{1 \le k \le \ell} \frac{a_\ell^k}{k!} + 2 \sum_{\ell < k \le N/2} b_\ell^k \\
 &\leq 2 \sum_{k \ge 1} \frac{a_\ell^k}{k!} + 2 \sum_{k \ge 1} b_\ell^k.
\end{align*}
Especially, when $b_\ell < 1$, we see that
\begin{equation}
 \label{eq:ExactRecoverysum}
 \err_z(0)
 \leq 2 (e^{a_\ell}-1) + \frac{2 b_\ell}{1-b_\ell}.
\end{equation}
Due to our choice $\ell = N^{2/3}$, we find that
\begin{align*}
 -\log a_\ell
 &= \frac12 (1-N^{-1/3}) N I_T - \log N, \\
 -\log b_\ell
 &= \frac14 N I_T - \frac13 \log N - 1.
\end{align*}
The assumption that $N I_T \ge (2+\delta) \log N$ for some scale-independent constant $\delta > 0$ now implies that $-\log a_\ell \to \infty$ and $-\log b_\ell \to \infty$, and therefore $a_\ell, b_\ell \to 0$.  Then \eqref{eq:ExactRecoverysum} shows that $\err_z(0) \to 0$ and confirms Theorem~\ref{the:ExactRecovery}.\hfill
\IEEEQED

\pagebreak
\bibliographystyle{IEEEtranN}
\bibliography{mybibliography}

\newcommand{\SortNoop}[1]{}\def\cprime{$'$}
\begin{thebibliography}{31}
\providecommand{\natexlab}[1]{#1}
\providecommand{\url}[1]{#1}
\csname url@samestyle\endcsname
\providecommand{\newblock}{\relax}
\providecommand{\bibinfo}[2]{#2}
\providecommand{\BIBentrySTDinterwordspacing}{\spaceskip=0pt\relax}
\providecommand{\BIBentryALTinterwordstretchfactor}{4}
\providecommand{\BIBentryALTinterwordspacing}{\spaceskip=\fontdimen2\font plus
\BIBentryALTinterwordstretchfactor\fontdimen3\font minus
  \fontdimen4\font\relax}
\providecommand{\BIBforeignlanguage}[2]{{%
\expandafter\ifx\csname l@#1\endcsname\relax
\typeout{** WARNING: IEEEtranN.bst: No hyphenation pattern has been}%
\typeout{** loaded for the language `#1'. Using the pattern for}%
\typeout{** the default language instead.}%
\else
\language=\csname l@#1\endcsname
\fi
#2}}
\providecommand{\BIBdecl}{\relax}
\BIBdecl

\bibitem[Holland et~al.(1983)Holland, Laskey, and
  Leinhardt]{holland_laskey_leinhardt_1983}
P.~W. Holland, K.~B. Laskey, and S.~Leinhardt, ``Stochastic blockmodels:
  {F}irst steps,'' \emph{Social Networks}, vol.~5, no.~2, pp. 109--137, 1983.

\bibitem[Hofman and Wiggins(2008)]{Hofman_Wiggins_2008}
J.~M. Hofman and C.~H. Wiggins, ``Bayesian approach to network modularity,''
  \emph{Physical Review Letters}, vol. 100, no.~25, p. 258701, 2008.

\bibitem[Peixoto(2019)]{Peixoto_2019}
T.~P. Peixoto, ``Bayesian stochastic blockmodeling,'' \emph{Advances in network
  clustering and blockmodeling}, pp. 289--332, 2019.

\bibitem[van~der Pas and van~der Vaart(2018)]{vanderpas_vandervaart_2018}
S.~van~der Pas and A.~van~der Vaart, ``Bayesian community detection,''
  \emph{Bayesian Analysis}, vol.~13, no.~3, pp. 767--796, 2018.

\bibitem[Ghosh et~al.(2019)Ghosh, Pati, and
  Bhattacharya]{ghosh_pati_bhattacharya_2019}
P.~Ghosh, D.~Pati, and A.~Bhattacharya, ``Posterior contraction rates for
  stochastic block models,'' \emph{Sankhya A}, vol.~82, no.~2, p. 448–476,
  2019.

\bibitem[Kleijn and van Waaij(2021)]{kleijn_van_waaij_2021}
B.~J.~K. Kleijn and J.~van Waaij, ``Confidence sets in a sparse stochastic
  block model with two communities of unknown sizes,'' 2021, arXiv:2108.07078.

\bibitem[Jiang and Tokdar(2021)]{jiang_tokdar_2021}
S.~Jiang and S.~Tokdar, ``Consistent {Bayesian} community detection,'' 2021,
  arXiv:2101.06531.

\bibitem[Li et~al.(2021)Li, Zhou, Tsui, Wei, and Ji]{li_zhou_tsui_wei_ji_2021}
T.~Li, T.~Zhou, K.-W. Tsui, L.~Wei, and Y.~Ji, ``Posterior contraction rate of
  sparse latent feature models with application to proteomics,''
  \emph{Statistical Theory and Related Fields}, pp. 1--11, 2021.

\bibitem[Herlau et~al.(2013)Herlau, M{\o}rup, and
  Schmidt]{herlau_morup_schmidt_2013}
T.~Herlau, M.~M{\o}rup, and M.~Schmidt, ``Modeling temporal evolution and
  multiscale structure in networks,'' in \emph{International Conference on
  Machine Learning}.\hskip 1em plus 0.5em minus 0.4em\relax PMLR, 2013, pp.
  960--968.

\bibitem[Xu and Hero(2014)]{xu_hero_2014}
K.~S. Xu and A.~O. Hero, ``Dynamic stochastic blockmodels for time-evolving
  social networks,'' \emph{IEEE Journal of Selected Topics in Signal
  Processing}, vol.~8, no.~4, pp. 552--562, 2014.

\bibitem[Han et~al.(2015)Han, Xu, and Airoldi]{han_xu_airoldi_2015}
Q.~Han, K.~Xu, and E.~Airoldi, ``Consistent estimation of dynamic and
  multi-layer block models,'' in \emph{Proceedings of the 32nd International
  Conference on Machine Learning}, ser. Proceedings of Machine Learning
  Research, vol.~37.\hskip 1em plus 0.5em minus 0.4em\relax PMLR, 2015, pp.
  1511--1520.

\bibitem[Matias and Miele(2017)]{matias_miele_2017}
C.~Matias and V.~Miele, ``Statistical clustering of temporal networks through a
  dynamic stochastic block model,'' \emph{Journal of the Royal Statistical
  Society B}, vol.~79, no.~4, pp. 1119--1141, 2017.

\bibitem[Ghasemian et~al.(2016)Ghasemian, Zhang, Clauset, Moore, and
  Peel]{ghasemian_zhang_clauset_moore_peel_2016}
A.~Ghasemian, P.~Zhang, A.~Clauset, C.~Moore, and L.~Peel, ``Detectability
  thresholds and optimal algorithms for community structure in dynamic
  networks,'' \emph{Physical Review X}, vol.~6, no.~3, p. 031005, 2016.

\bibitem[Barucca et~al.(2018)Barucca, Lillo, Mazzarisi, and
  Tantari]{Barucca_Lillo_Mazzarisi_Tantari_2018}
P.~Barucca, F.~Lillo, P.~Mazzarisi, and D.~Tantari, ``Disentangling group and
  link persistence in dynamic stochastic block models,'' \emph{Journal of
  Statistical Mechanics: Theory and Experiment}, vol. 2018, no. 123407, pp.
  1--18, 2018.

\bibitem[Bhattacharyya and
  Chatterjee(2020)]{Bhattacharyya_Chatterjee_2020-04-06}
S.~Bhattacharyya and S.~Chatterjee, ``General community detection with optimal
  recovery conditions for multi-relational sparse networks with dependent
  layers,'' 2020, arXiv:2004.03480.

\bibitem[Durante et~al.(2017)Durante, Dunson, and
  Vogelstein]{durante_dunson_vogelstein_2017}
D.~Durante, D.~B. Dunson, and J.~T. Vogelstein, ``Nonparametric bayes modeling
  of populations of networks,'' \emph{Journal of the American Statistical
  Association}, vol. 112, no. 520, pp. 1516--1530, 2017.

\bibitem[Lei et~al.(2020)Lei, Chen, and Lynch]{lei_chen_lynch_2020}
J.~Lei, K.~Chen, and B.~Lynch, ``Consistent community detection in multi-layer
  network data,'' \emph{Biometrika}, vol. 107, no.~1, pp. 61--73, 2020.

\bibitem[Avrachenkov et~al.(2022)Avrachenkov, Dreveton, and
  Leskel\"a]{Avrachenkov_Dreveton_Leskela_2022}
K.~Avrachenkov, M.~Dreveton, and L.~Leskel\"a, ``Community recovery in
  non-binary and temporal stochastic block models,'' 2022, arXiv:2008.04790.

\bibitem[{van Erven} and {Harremo\"es}(2014)]{vanErven_Harremoes_2014}
T.~{van Erven} and P.~{Harremo\"es}, ``R{\'e}nyi divergence and
  {Kullback}--{Leibler} divergence,'' \emph{IEEE Transactions on Information
  Theory}, vol.~60, no.~7, pp. 3797--3820, 2014.

\bibitem[Abbe et~al.(2015)Abbe, Bandeira, and Hall]{abbe_bandeira_hall_2015}
E.~Abbe, A.~S. Bandeira, and G.~Hall, ``Exact recovery in the stochastic block
  model,'' \emph{IEEE Transactions on information theory}, vol.~62, no.~1, pp.
  471--487, 2015.

\bibitem[Mossel et~al.(2015)Mossel, Neeman, and Sly]{Mossel_Neeman_Sly_2016}
E.~Mossel, J.~Neeman, and A.~Sly, ``Consistency thresholds for the planted
  bisection model,'' in \emph{Proceedings of the 47th Annual ACM symposium on
  Theory of Computing}, 2015, pp. 69--75.

\bibitem[Nowicki and Snijders(2001)]{nowicki_snijders_2001}
K.~Nowicki and T.~A.~B. Snijders, ``Estimation and prediction for stochastic
  blockstructures,'' \emph{Journal of the American Statistical Association},
  vol.~96, no. 455, p. 1077–1087, 2001.

\bibitem[Gao et~al.(2015)Gao, Lu, and Zhou]{gao_lu_zhou_2015}
C.~Gao, Y.~Lu, and H.~H. Zhou, ``Rate-optimal graphon estimation,''
  \emph{Annals of Statistics}, vol.~43, no.~6, pp. 2624--2652, 2015.

\bibitem[Klopp et~al.(2017)Klopp, Tsybakov, and
  Verzelen]{klopp_tsybakov_verzelen_2017}
O.~Klopp, A.~B. Tsybakov, and N.~Verzelen, ``Oracle inequalities for network
  models and sparse graphon estimation,'' \emph{Annals of Statistics}, vol.~45,
  no.~1, pp. 316--354, 2017.

\bibitem[Dahl(2006)]{dahl_2006}
D.~B. Dahl, ``Model-based clustering for expression data via a dirichlet
  process mixture model,'' \emph{Bayesian inference for gene expression and
  proteomics}, vol.~4, pp. 201--218, 2006.

\bibitem[Rand(1971)]{rand_1971}
W.~M. Rand, ``Objective criteria for the evaluation of clustering methods,''
  \emph{Journal of the American Statistical Association}, vol.~66, no. 336, pp.
  846--850, 1971.

\bibitem[Hubert and Arabie(1985)]{hubert_arabie_1985}
L.~Hubert and P.~Arabie, ``Comparing partitions,'' \emph{Journal of
  Classification}, vol.~2, no.~1, pp. 193--218, 1985.

\bibitem[G{\"o}sgens et~al.(2021)G{\"o}sgens, Tikhonov, and
  Prokhorenkova]{Gosgens_Tikhonov_Prokhorenkova_2021}
M.~M. G{\"o}sgens, A.~Tikhonov, and L.~Prokhorenkova, ``Systematic analysis of
  cluster similarity indices: How to validate validation measures,'' in
  \emph{Proceedings of the 38th International Conference on Machine Learning},
  vol. 139.\hskip 1em plus 0.5em minus 0.4em\relax PMLR, 18--24 Jul 2021, pp.
  3799--3808.

\bibitem[Lei et~al.(2017)Lei, Bezdek, Romano, Vinh, Chan, and
  Bailey]{Lei_etal_2017}
Y.~Lei, J.~C. Bezdek, S.~Romano, N.~X. Vinh, J.~Chan, and J.~Bailey, ``Ground
  truth bias in external cluster validity indices,'' \emph{Pattern
  Recognition}, vol.~65, pp. 58--70, 2017.

\bibitem[Ghosal and van~der Vaart(2017)]{Ghosal_VanDerVaart_2017}
S.~Ghosal and A.~van~der Vaart, \emph{Fundamentals of nonparametric Bayesian
  inference}.\hskip 1em plus 0.5em minus 0.4em\relax Cambridge University
  Press, 2017.

\bibitem[Kleijn(2021)]{Kleijn_2021}
B.~Kleijn, ``Frequentist validity of {Bayesian} limits,'' \emph{Annals of
  Statistics}, vol.~49, no.~1, pp. 182--202, 2021.

\end{thebibliography}

\end{document}